%% file: spcavity.tex
\documentclass[10pt,a4paper,oneside]{article}

\usepackage{xcolor}

\usepackage{indentfirst}

\usepackage{mathtools} 
\usepackage{amssymb} 
\usepackage{amsthm}  
\usepackage{stmaryrd} 
\usepackage{accents} 

\usepackage[normalem]{ulem} 
\usepackage{soul} 

\usepackage{array}
\usepackage{arydshln} 
\usepackage{booktabs} 
\usepackage{multirow} 
\usepackage{enumitem} 

\usepackage{graphicx}
\usepackage[labelformat=simple]{subcaption}
\usepackage{tikz}

\usepackage[algo2e,vlined,ruled,algosection,linesnumberedhidden]{algorithm2e}

\usepackage{fancyhdr}
\usepackage{lastpage}

\usepackage[colorlinks, linkcolor=blue, citecolor=blue]{hyperref}

\newtheorem{theorem}{Theorem}[section]
\newtheorem{lemma}[theorem]{Lemma}

\numberwithin{equation}{section}
\numberwithin{table}{section}
\numberwithin{figure}{section}

\input{macros}

\ProvideDocumentEnvironment{keywords}{}
  {\par\noindent\textbf{Keywords: }}
  {\par}
 
\begin{document}

\title{Scattering Resonances of Plane-Embedded Cavities: A Finite Element DtN Approach}

\author{Bo Gong$^1$, Takumi Sato$^2$, Jiguang Sun$^3$, Xinming Wu$^{2,}$\thanks{Corresponding author. Email address: wuxinming@fudan.edu.cn}}

\footnotetext[1]{School of Mathematics, Statistics and Mechanics, Beijing University of Technology, Beijing, 100124, China.}
\footnotetext[2]{School of Mathematical Sciences, SKLCAM, Fudan University, Shanghai 200433, China.}
\footnotetext[3]{Department of Mathematical Sciences, Michigan Technological University, Houghton, MI 49931, U.S.A.}

\date{}

\maketitle

\begin{abstract}
While scattering resonances of plane-embedded cavities have wide-ranging engineering applications, numerical results for such configurations remain virtually unexplored. We present an effective finite element approach for computing these resonances. The unbounded domain is truncated using a Dirichlet-to-Neumann (DtN) boundary operator, recasting the problem as the eigenvalue problem of a holomorphic Fredholm operator function. Discretization is carried out using Lagrange finite elements with a truncated DtN expansion, and convergence is rigorously proven using Karma’s approximation theory. The resulting non-linear matrix eigenvalue system is computed using a multi-step parallel spectral indicator method. Numerical experiments confirm the theoretical predictions and provide benchmark solutions for future investigations. 
\end{abstract}


\section{Introduction}
Scattering resonances of an open cavity are complex frequencies at which the cavity supports nontrivial outgoing waves without an incident field. They are also called  scattering poles or quasinormal-mode frequencies \cite{LaxPhillips1989, Zworski1999, DyatlovZworski2019}. Suppose a wave enters an open cavity through its aperture. It undergoes repeated reflections from the cavity walls, and part of the wave escapes after each reflection. At certain frequencies, the internally reflected waves interfere constructively and produce a strong localized field. Because the cavity is open, energy leaks into the exterior. 

Scattering resonances of cavities play an important role in predicting and controlling wave interactions with open or recessed structures. They are used in the design of cavity-backed antennas, acoustic noise control in ducts and mufflers, and the development of sensors and resonators. They are also useful for detecting cracks, voids, corrosion, tunnels, and other hidden defects through inverse scattering and nondestructive testing \cite{Wood06, Bao2014, Ammari2020, LinZhang2026}.

Despite of extensive analysis, there seem to be no rigorous numerical analysis results on the scattering resonances of open cavities. It is the goal of this paper to propose a finite element method and prove the error estimates. The unbounded domain is first truncated and a Dirichlet-to-Neumann (DtN) mapping  is applied on the artificial boundary. The DtN mapping is truncated and problem is discretized using Lagrange elements. The resulting nonlinear matrix eigenvalue problem is computed using the parallel multi-step spectral indicator method.  The main tool in the error analysis is the abstract approximation theory for eigenvalue problems of holomorphic Fredholm operator functions \cite{Karma1, Karma2}. 

The framework of the combination of DtN mapping, finite elements, and spectral indicator method has been used to successfully compute the scattering resonances for several problems including the acoustic and electromagnetic scattering of obstacles and inhomogeneous media \cite{XiGongSun2024, Gong2026AML, XiLinSun2024CMA, Gong2026a, Gong2026b}. The method does not produce non-physical (spurious) resonances. In contrast, the use of PML computes a non-physical essential spectrum and other spurious resonances \cite{KimPasciak2009MC, Araujo2021JCP}. Boundary integral equations are usually not equivalent to the original problem and produce non-physical resonances, in particular, in the case of inhomogeneous media \cite{NannenWess2018BIT, Grubisic2023}. 

The contributions of the paper include: 1) a spectrum correct computational method for scattering resonances of a cavity embedded in the infinite plane, which integrates finite elements, DtN mapping, and the spectral indicator methods; 2) error estimates taking both the truncation error of DtN series and finite element error into account; 3) benchmarks for both overfilled and empty cavities. The rest of the paper is organized as follows. In section~\ref{Sec2}, we present the scattering problem defined on the unbounded domain and an equivalent problem on the truncated domain terminated by the DtN mapping. We show that the resonances are the eigenvalues of a holomorphic Fredholm operator function associated to the weak formulation of the truncated problem. In Section~\ref{FEM}, we propose a finite element method and prove the error estimates. Section~\ref{NE} is devoted to extensive numerical study of non-filled or over-filled cavities. 

\section{Scattering resonances of cavities}\label{Sec2}

\begin{figure}[htbp!]
  \centering
  \includegraphics[width=0.55\textwidth]{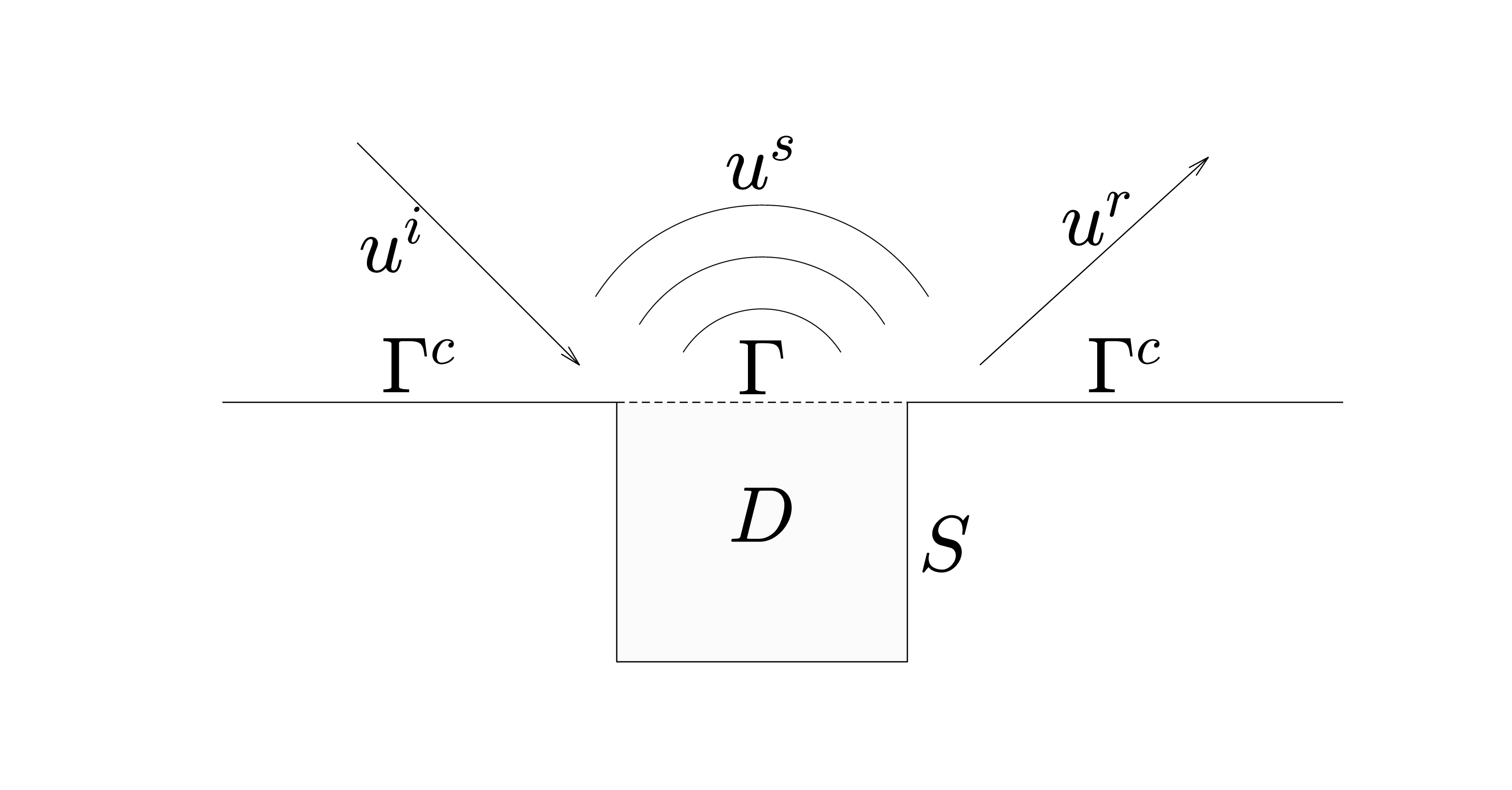}
  \caption{Schematic of the scattering problem for an open cavity.}
  \label{fig:cavity-geometry}
\end{figure}

Let $\bbR^2_+$ and $\bbR^2_-$ be the upper and lower half-planes, respectively. Let $D\subset \bbR^2_-$ denote a cavity, a bounded Lipschitz domain, embedded in the infinite half plane. The boundary of $D$ consists of the open aperture $\Gamma$, a line segment, and the wall $S$. Let $\Gamma^c=\partial\bbR^2_+\setminus\overline\Gamma$ (see Fig.~\ref{fig:cavity-geometry}). Consider the time-harmonic scattering problem of a cavity 
\begin{subequations}
\begin{align}
    &\Delta u + k^2n(x)u = 0 &&\text{in } D\cup\bbR^2_+, \label{eq:sc1}\\
    &u = 0 &&\text{on } \Gamma^c\cup S, \label{eq:sc2}\\
    &\sqrt{r}\left(\frac{\partial u^s}{\partial r} - \rmi ku^s\right) \to 0 
    &&\text{as } r=|x|\to \infty, \label{eq:sc3}
\end{align}
\end{subequations}
where $k$ is the wavenumber and $n(x)$ is the index of refraction.

When $n(x)=1$, the cavity is non-filled. The cavity can also be filled or overfilled, meaning that the medium fills the cavity or even extends above the aperture $\Gamma$ into $\bbR^2_+$.
For an incident field given by
\begin{equation}
    u^{i}=e^{\rmi k(x\sin\theta-y\cos\theta)},
\end{equation}
where $\theta\in(-\pi/2,\pi/2)$ denotes the angle of incidence measured from the negative $y$-axis, the corresponding reflected field for the half-plane is given by
\begin{equation}
    u^{r}=-e^{\rmi k(x\sin\theta+y\cos\theta)}.
\end{equation}
The total field
\begin{equation}
    u=u^{i}+u^{r}+u^s
\end{equation}
satisfies the Helmholtz equation \eqref{eq:sc1} and the homogeneous Dirichlet boundary condition \eqref{eq:sc2}, which corresponds to the perfectly conducting boundary condition. The scattered field $u^s$ satisfies the Sommerfeld radiation condition \eqref{eq:sc3} in the upper half-plane $\bbR^2_+$. 

If the imaginary part of $k$ is non-negative, i.e., $\Im(k)\ge 0$, the scattering problem \eqref{eq:sc1}-\eqref{eq:sc3} has a unique solution $u\in H^1_{\rm loc}(D\cup\bbR^2_+)$ (see, e.g., \cite{LW13, Wood06}). The solution operator function $S(k): u^i\to u$ is holomorphic in $k$ on the upper half complex plane and can be meromorphically extended to $\bbC\bs\bbR_-$ with $\bbR_-$ being the set of non-positive real numbers. The scattering resonances are the poles of $S(k)$ in $\{k\in\bbC : \Im(k) < 0\}$ \cite{LaxPhillips1989, DyatlovZworski2019}..

We now introduce an outgoing condition for the scattered field in the upper half-plane $\bbR^2_+$. Suppose that the aperture $\Gamma$ is contained in the boundary segment $\{(x,0): x\in\bbR, \,|x|<R\}$. For $r=|x|>R$, the homogeneous Dirichlet boundary conditions at $\theta=0$ and $\theta=\pi$ define an angular eigenvalue problem whose eigenfunctions are $\sin(n\theta),\ n\geq 1$. These eigenfunctions form a complete orthogonal basis in $L^2(0,\pi)$. Expanding the scattered field in this basis yields the outgoing wave expansion
\begin{equation}
    \label{eq:outgoing}
    u^s(r,\theta)=\sum_{n=1}^\infty a_n H_n^{(1)}(kr)\sin n\theta, 
    \qquad r>R,\quad 0<\theta<\pi,
\end{equation}
where $H_n^{(1)}$ denotes the Hankel function of the first kind of order $n$.
Moreover, for $k$ with $\Im(k)\ge 0$, \eqref{eq:outgoing} is equivalent to the Sommerfeld radiation condition \eqref{eq:sc3} together with the condition $u^s=0$ on $\Gamma^c$. It remains valid for complex wavenumbers $k$ with $\Im(k)<0$.

\subsection{Weak formulation and operator equation}

To reduce the scattering problem \eqref{eq:sc1}-\eqref{eq:sc3} to a bounded domain, we introduce the following Dirichlet-to-Neumann (DtN) map. Let $R>0$ be fixed and
\begin{equation*}
    B_R^+=\{(r,\theta): 0\le r<R,\ \theta\in(0,\pi)\}
\end{equation*}
be the upper half-disk of radius $R$ and $\Gamma_R^+$ denote its semicircular boundary, i.e., $\Gamma_R^+=\{(r,\theta): r=R,\ \theta\in(0,\pi)\}$. For $s\ge 0$, the Sobolev space associated with the Fourier sine basis is defined as
\begin{equation*}
    H^s_D(\Gamma_R^+)
    := \left\{ \varphi=\sum_{n=1}^{\infty}\varphi_n\sin(n\theta):
    \sum_{n=1}^{\infty}(1+n^2)^s|\varphi_n|^2<\infty \right\},
\end{equation*}
equipped with the norm
\begin{equation}
    \|\varphi\|_{H_D^s(\Gamma_R^+)}^2
    := \frac{\pi R}{2} \sum_{n=1}^{\infty}(1+n^2)^s|\varphi_n|^2.
\end{equation}
The negative-order space is defined by duality as
\begin{equation*}
    H_D^{-s}(\Gamma_R^+)
    :=\bigl(H_D^s(\Gamma_R^+)\bigr)', \quad s>0
\end{equation*}
with $L^2(\Gamma_R^+)$ as the pivot space. Equivalently,
\begin{equation*}
    \|\varphi\|_{H_D^{-s}(\Gamma_R^+)}^2
    = \frac{\pi R}{2} \sum_{n=1}^{\infty}(1+n^2)^{-s}|\varphi_n|^2,
\end{equation*}
where the Fourier sine expansion is understood in the distributional sense.


By the orthogonality of $\{\sin n\theta\}_{n=1}^\infty$ on $(0,\pi)$, the coefficients in \eqref{eq:outgoing} are given by
\begin{equation}
    a_n=\frac{2}{\pi H_n^{(1)}(kR)} 
    \int_0^\pi u^s(R,\theta')\sin n\theta'\,\rmd\theta'.
\end{equation}
Hence, for $r\ge R$,
\begin{equation}
    u^s(r,\theta) = \frac{2}{\pi} \sum_{n=1}^\infty \frac{H_n^{(1)}(kr)}{H_n^{(1)}(kR)} \sin n\theta 
    \int_0^\pi u^s(R,\theta')\sin n\theta'\,\rmd\theta'.
\end{equation}
Consequently, the DtN map $T(k): H_D^{1/2}(\Gamma_R^+)\to H_D^{-1/2}(\Gamma_R^+)$ is defined as
\begin{equation}
    \label{eq:operatorT}
    T(k)\varphi = \frac{2k}{\pi} \sum_{n=1}^\infty \frac{H_n^{(1)'}(kR)}{H_n^{(1)}(kR)} \sin n\theta 
    \int_0^\pi \varphi(\theta')\sin(n\theta')\,\rmd\theta',
\end{equation}
where $\varphi(\theta)=u^s(R,\theta)$ is the Dirichlet trace of $u^s$ on $\Gamma_R^+$.

Let $\calZ_n$ be the zeros of $H_n^{(1)}(kR)$ and $\calZ = \cup_{n=1}^\infty \calZ_n$. Define $\Lambda=\bbC\bs(\bbR_-\cup\calZ)$, for $k\in\Lambda$, the DtN map $T(k)$ defined in \eqref{eq:operatorT} is bounded, i.e.
\begin{equation}
    \label{eq:Tbound}
    \|T(k)\varphi\|_{H_D^{-1/2}(\Gamma_R^+)} \le C \|\varphi\|_{H_D^{1/2}(\Gamma_R^+)},
\end{equation}
where $C$ is a constant depending on $k$ but not on $\varphi$ (\cite{Wood06}). 
Using $T(k)$, the scattering problem \eqref{eq:sc1}--\eqref{eq:sc3} can be formulated on the bounded domain $\Omega=D\cup B_R^+$. Let
\begin{equation*}
    \Gamma_R^{\,c}:=\{(x,0): x\in\bbR,\,|x|<R\}\setminus\overline{\Gamma}
\end{equation*}
be the truncated part of the perfectly conducting boundary outside the aperture. Then the truncated problem is given by
\begin{subequations}
\begin{align}
    &\Delta u + k^2n(x)u = 0 &&\text{in } \Omega, \label{eq:u1}\\
    &u = 0 &&\text{on } \Gamma_R^{\,c}\cup S, \label{eq:u2}\\
    &\frac{\partial u}{\partial \nu} - T(k)u = g(k) &&\text{on } \Gamma_R^+, \label{eq:u3}
\end{align}
\end{subequations}
where $\nu$ is the unit outward normal to $\Gamma_R^+$ and $g(k)=\frac{\partial u^i}{\partial \nu} + \frac{\partial u^r}{\partial \nu}  - T(k)(u^i+u^r)$.

Define the space
\begin{equation}
    V =\{\varphi\in H^1(\Omega)\,:\, \varphi=0\;\text{on}\ \Gamma_R^{\,c}\cup S\}.
\end{equation}
Given $g(k)\in H_D^{-1/2}(\Gamma_R^+)$, the variational formulation for the truncated problem \eqref{eq:u1}-\eqref{eq:u3} is to find $u\in V$ such that
\begin{equation}
\label{eq:weak}
    (\nabla u,\nabla v) - k^2(n(x)u,v) - \langle T(k)u, v\rangle_{\Gamma_R^+} 
    = \langle g(k), v\rangle_{\Gamma_R^+} \quad\forall v\in V,
\end{equation}
where $(\cdot,\cdot)$ denotes the $L^2(\Omega)$ inner product, $\langle\cdot,\cdot\rangle_{\Gamma_R^+}$ denotes the $H_D^{-1/2}(\Gamma_R^+)$-$H_D^{1/2}(\Gamma_R^+)$ duality.

For any $k\in\Lambda$, define the operator $B(k):V\to V$ such that
\begin{equation}
    \label{eq:operatorB}
    (B(k)u, v)_{H^1(\Omega)} = (\nabla u,\nabla v) - k^2(n(x)u,v) 
    - \langle T(k)u, v\rangle_{\Gamma_R^+} \quad\forall u,v\in V,
\end{equation}
where $(\cdot,\cdot)_{H^1(\Omega)}$ is the $H^1(\Omega)$ inner product.
Define $f(k)\in V$ such that
\begin{equation}
    (f(k),v)_{H^1(\Omega)} = \langle g(k),v\rangle_{\Gamma_R^+} \quad \forall v\in V.
\end{equation}
Then \eqref{eq:weak} can be written equivalently as
\begin{equation}
    B(k)u=f(k) \quad \text{in } V.
\end{equation}
The scattering poles are the values $k\in\Lambda$ at which the meromorphic operator-valued function $B(\cdot)^{-1}$ has a pole. Equivalently, they are
those $k\in\Lambda$ for which $B(k)$ is not invertible. Hence, the associated operator eigenvalue problem is to find $(k,u)\in\Lambda\times V$ such that
\begin{equation}
    \label{eq:eig}
    B(k)u=0 \quad \text{in } V.
\end{equation}

\subsection{Holomorphic Fredholm operator functions}
The operator $B(k)$ is a holomorphic Fredholm operator function, for which we recall some basics (see e.g. \cite{Karma1}). Let $X$ and $Y$ be Banach spaces. An operator $F\in \mathcal{L}(X,Y)$ is Fredholm if the dimension of the nullspace $\mathcal{N}(F)$ is finite, the range $\mathcal{R}(F)$ is closed, and the codimension of $\mathcal{R}(F)$ is finite. The index of a Fredholm operator $F$ is defined as $\text{dim}\mathcal{N}(F) - \text{codim}\mathcal{R}(F)$. We say that an operator function $F(\cdot) : \Lambda\to \mathcal{L}(X,Y)$ is a Fredholm operator function (with index zero), if $F(\lambda)$ is Fredholm (with index zero) for all $\lambda\in \Lambda$.

A tuple $(x_0, x_1, \dots, x_j)$ with $x_0\neq 0$ and
\begin{equation*}
    F(\lambda)x_\ell + F^{(1)}(\lambda)x_{\ell-1}
    + \cdots +\frac{1}{\ell!} F^{(\ell)}(\lambda)x_0 = 0,
    \quad \ell=0,1,\dots,j
\end{equation*}
is called a Jordan chain of length $j+1$ of $F(\cdot)$ at $\lambda$. Such elements $x_0, x_1, \dots, x_j$ are called
the generalized eigenelements of $F(\cdot)$ at $\lambda$,
and the closed span of all the generalized eigenelements
of $F(\cdot)$ at $\lambda$ is called the generalized eigenspace
of $F(\cdot)$ at $\lambda$, denoted as $G(F,\lambda)$.

If $F(\cdot):\Lambda\to \mathcal{L}(X,Y)$ is a Fredholm operator function with index zero, and if the resolvent set
$\rho(F) = \{\lambda\in \Lambda\,:\, F(\lambda)^{-1}\in \mathcal{L}(Y,X)\}$ is non-empty, then the spectrum $\sigma(F) = \lambda\setminus \rho(F)$ of $F(\cdot)$ has no cluster point in $\Lambda$ and every $\lambda\in \sigma(F)$ is an eigenvalue of $F(\cdot)$, furthermore, the generalized eigenspace $G(F,\lambda)$ is finite dimensional.

\subsection{Holomorphy and Fredholmness of $B(\cdot)$}

We now analyze the operator function $B(\cdot)$.
For any $u\in V$, we write their traces on $\Gamma_R^+$ in the sine basis as
\begin{equation}
    u(R,\theta)=\sum_{n=1}^{\infty}u_n\sin n\theta,
\end{equation}
where
\begin{equation}
    u_n=\frac{2}{\pi}\int_0^\pi u(R,\theta)\sin n\theta\,d\theta.
\end{equation}
The DtN map $T(k)$ in \eqref{eq:operatorT} can be written as
\begin{equation}
    \label{eq:operatorT2}
    T(k) u = \sum_{n=1}^{\infty} \tau_n(k) u_n \sin n\theta,\quad \text{where}\quad \tau_n(k) = k\frac{H_n^{(1)'}(kR)}{H_n^{(1)}(kR)}.
\end{equation}
The duality product is given by
\begin{equation}
    \label{eq:Tuv}
    \langle T(k)u,v\rangle_{\Gamma_R^+}
    = \frac{\pi R}{2} \sum_{n=1}^{\infty} \tau_n(k) u_n\overline{v_n}.
\end{equation}
Using the recurrence relation
\begin{equation}
    \label{eq:hprime}
    H_n^{(1)'}(z) = H_{n-1}^{(1)}(z)- \frac{n}{z}H_n^{(1)}(z),
\end{equation}
we obtain that
\begin{equation}
    \label{eq:tau}
    \tau_n(k) = \ -\frac{n}{R} + k\frac{H_{n-1}^{(1)}(kR)}{H_n^{(1)}(kR)}.
\end{equation}

Consequently, for each $k\in\Lambda$, the operator $B(k)$ can be decomposed as
\begin{equation}
    B(k)=A+K(k),
\end{equation}
where $A:V\to V$ and $K(k):V\to V$ are respectively defined via the $H^1(\Omega)$ inner product by
\begin{equation}\label{eq:operatorA}
    (Au,v)_{H^1(\Omega)}
    = (\nabla u,\nabla v)+(u,v) 
    + \frac{\pi}{2}\sum_{n=1}^{\infty} nu_n\overline{v_n} 
\end{equation}
and
\begin{equation}\label{eq:operatorK}
    (K(k)u,v)_{H^1(\Omega)}
    =-((k^2n(x)+1)u,v) 
    - \frac{\pi}{2}\sum_{n=1}^{\infty} kR\frac{H_{n-1}^{(1)}(kR)}{H_n^{(1)}(kR)} u_n\overline{v_n}
\end{equation}
for all $u,v\in V$.
%
%

\begin{theorem}\label{thm:B}
The operator function $B(\cdot): \Lambda\to\calL(V)$ is holomorphic and Fredholm with index zero. Furthermore, the resolvent set $\rho(B)$ is non-empty.
\end{theorem}

\begin{proof}
The operator $A$ is bounded and coercive on $V$. Indeed, by the trace theorem
and the characterization of the $H_D^{1/2}$ norm on $\Gamma_R^+$ in terms of
the sine coefficients, the boundary term in \eqref{eq:operatorA} is bounded on
$V\times V$. Moreover,
\begin{equation*}
    (Au,u)_{H^1(\Omega)}
    = \|\nabla u\|_{L^2(\Omega)}^2 + \|u\|_{L^2(\Omega)}^2
    + \frac{\pi}{2} \sum_{n=1}^{\infty} n |u_n|^2
    \ge \|u\|_{H^1(\Omega)}^2 .
\end{equation*}
Hence, by the Lax--Milgram theorem, $A:V\to V$ is an isomorphism.

We next show that $K(k):V\to V$ is compact for each $k\in\Lambda$. The volume integral part of $K(k)$ is compact because the embedding $V\hookrightarrow L^2(\Omega)$ is compact. For the boundary part, for each fixed $k\in\Lambda$, the large-order asymptotic behavior of Hankel functions gives
\begin{equation}
    \label{eq:hdh}
    \frac{H_{n-1}^{(1)}(kR)}{H_n^{(1)}(kR)}
    = O\left(\frac{1}{n}\right),\quad n\to\infty.
\end{equation}
After enlarging the constant to absorb finitely many lower-order modes, there exists a constant $C_k>0$ such that
\begin{equation}
    \label{eq:hdh2}
    \left| \frac{H_{n-1}^{(1)}(kR)}{H_n^{(1)}(kR)} \right|
    \le \frac{C_k}{n},\quad n\ge 1.
\end{equation}
Therefore, using the characterization of Sobolev norms on $\Gamma_R^+$ in
terms of sine coefficients, we obtain
\begin{equation*}
    \left| \sum_{n=1}^{\infty} kR \frac{H_{n-1}^{(1)}(kR)}{H_n^{(1)}(kR)} u_n\overline{v_n} \right|
    \le C_k \|u\|_{H_D^{-3/2}(\Gamma_R^+)} \|v\|_{H_D^{1/2}(\Gamma_R^+)}.
\end{equation*}
Define
\[
H_{00}^{1/2}(\Gamma_R^+)
:=
\left\{
\varphi\in H^{1/2}(\Gamma_R^+):
\widetilde{\varphi}\in H^{1/2}(\partial\Omega)
\right\},
\quad
\|\varphi\|_{H_{00}^{1/2}(\Gamma_R^+)}
:=
\|\widetilde{\varphi}\|_{H^{1/2}(\partial\Omega)},
\]
where \(\widetilde{\varphi}\) denotes the extension of \(\varphi\) by zero to \(\Gamma^c_R\cup S\). By the trace theorem and the homogeneous boundary condition on \(\Gamma^c_R\cup S\), the trace operator from \(V\) into \(H_{00}^{1/2}(\Gamma_R^+)\) is bounded. Moreover, \(H_D^{1/2}(\Gamma_R^+)=H_{00}^{1/2}(\Gamma_R^+)\), with equivalent norms, hence the trace operator from \(V\) into \(H_D^{1/2}(\Gamma_R^+)\) is also bounded.
Meanwhile, the embedding $H_D^{1/2}(\Gamma_R^+)\hookrightarrow H_D^{-3/2}(\Gamma_R^+)$ is compact.
Therefore the boundary part of $K(k)$ defines a compact operator on $V$. Consequently,
$B(k)=A+K(k)$ is a compact perturbation of the isomorphism $A$, i.e.,
$B(k)$ is Fredholm with index zero for every $k\in\Lambda$. Since $B(k)$ is
invertible for $k$ in the upper half-plane, its resolvent set $\rho(B)$ is
non-empty.

It remains to justify the holomorphic dependence on $k$. For any $u,v\in H_D^{1/2}(\Gamma_R^+)$ and any compact subset $\Lambda_0\Subset\Lambda$, by \eqref{eq:tau} and \eqref{eq:hdh2}, there exists a constant $C_{\Lambda_0}>0$ such that
\begin{equation}
    \label{eq:tau2}
    \left| \tau_n(k) \right| \le C_{\Lambda_0}n
    \qquad \forall k\in\Lambda_0, \ n\ge 1.
\end{equation}
Since
\begin{equation}
    \sum_{n=1}^\infty n|u_n||v_n|
    \le \|u\|_{H_D^{1/2}(\Gamma_R^+)} \|v\|_{H_D^{1/2}(\Gamma_R^+)}
\end{equation}
is a convergent series, the Weierstrass $M$-test implies that the series in \eqref{eq:Tuv} converges uniformly on $\Lambda_0$. Since each term is holomorphic in $k$ and $\Lambda_0\Subset\Lambda$ is arbitrary, the Weierstrass theorem implies that the function $k\mapsto \langle T(k)u,v\rangle_{\Gamma_R^+}$ is holomorphic on $\Lambda$. Hence by \cite[Theorem 3.12]{Kato}, $T(\cdot)$ is holomorphic on $\Lambda$.
Consequently, by the definition of $B(k)$ in \eqref{eq:operatorB}, we conclude that $B(\cdot):\Lambda\to\mathcal L(V)$ is holomorphic.
\end{proof}

\subsection{Truncated DtN map}

For numerical implementation, the infinite series in \eqref{eq:operatorT2} for $T(k)$ needs be truncated. This yields a finite-term approximation
\begin{equation}
    \label{eq:operatorTN}
    T^N(k)u = \sum_{n=1}^N \tau_n(k) u_n \sin n\theta,
\end{equation}
where $N$ is referred to as the truncation order. Similar to \eqref{eq:Tbound}, the operator $T^N(k)$ is also bounded, i.e.,
\begin{equation}
    \label{eq:TNbound}
    \|T^N(k)\varphi\|_{H_D^{-1/2}(\Gamma_R^+)} \le C \|\varphi\|_{H_D^{1/2}(\Gamma_R^+)}.
\end{equation}

\begin{lemma}\label{lem:Truncation}
For $T(k)$ and $T^N(k)$ defined by \eqref{eq:operatorT} and \eqref{eq:operatorTN}, respectively, it holds that
\begin{equation}
    \label{eq:errT}
    \|(T(k)-T^N(k))\varphi\|_{H_D^{-1/2}(\Gamma_R^+)} 
    \le C\frac{\gamma_{1/2+s}(N,\varphi)}{N^s} \|\varphi\|_{H_D^{1/2+s}(\Gamma_R^+)}
\end{equation}
for all $\varphi\in H_D^{1/2+s}(\Gamma_R^+)$ with $s\ge0$, where $C > 0$ is a constant depending on $k$ but independent of $\varphi$ and $N$. Furthermore, $\gamma_{1/2+s}(N,\varphi) \in [0,1]$ and $\gamma_{1/2+s}(N,\varphi)\to 0$ as $N\to\infty$.
\end{lemma}

\begin{proof}
By the definitions \eqref{eq:operatorT2} and \eqref{eq:operatorTN},
\begin{equation*}
    (T(k)-T^N(k))\varphi = \sum_{n=N+1}^{\infty}\tau_n(k)\varphi_n\sin n\theta.
\end{equation*}
Using the estimate \eqref{eq:tau2}, we have that
\begin{align*}
    &\|(T(k)-T^N(k))\varphi\|_{H_D^{-1/2}(\Gamma_R^+)}^2
    = \sum_{n=N+1}^{\infty}(1+n^2)^{-1/2}|\tau_n(k)|^2|\varphi_n|^2\\
   \le{} & C\sum_{n=N+1}^{\infty}(1+n^2)^{1/2}|\varphi_n|^2
    \le \frac{C}{N^{2s}}\sum_{n=N+1}^{\infty}(1+n^2)^{1/2+s}|\varphi_n|^2 \\
    \le{}& \frac{C}{N^{2s}}\gamma_{1/2+s}^2(N,\varphi)\|\varphi\|_{H_D^{1/2+s}(\Gamma_R^+)}^2,
\end{align*}
where 
\[
\gamma_{1/2+s}(N,\varphi) = \Big\{\sum_{n=N+1}^\infty (1+n^2)^{1/2+s} |\varphi_n|^2\left/\sum_{n=1}^\infty (1+n^2)^{1/2+s} |\varphi_n|^2\right. \Big\}^{1/2}.
\]
\end{proof}

Using the truncated DtN operator $T^N(k)$, the approximation $B^N(k)$ of $B(k)$ is defined by
\begin{equation}
    \label{eq:operatorBN}
    (B^N(k)u, v)_{H^1(\Omega)} = (\nabla u,\nabla v) - k^2(n(x)u,v) 
    - \langle T^N(k)u, v\rangle_{\Gamma_R^+} \quad\forall  u,v\in V.
\end{equation}
The associated eigenvalue problem is to find $k$ and $u^N$ such that
\begin{equation}
    \label{eq:eigN}
    B^N(k)u^N=0 \quad \text{in } V.
\end{equation}
The operator $B^N(k)$ can be split analogously
\begin{equation}
    B^N(k)=A+K^N(k),
\end{equation}
where $A^N:V\to V$ and $K^N(k):V\to V$ are respectively defined via the $H^1(\Omega)$-inner product by
\begin{equation}
    (A^Nu,v)_{H^1(\Omega)}
    = (\nabla u,\nabla v)+(u,v) 
    + \frac{\pi}{2}\sum_{n=1}^{N} nu_n\overline{v_n} \quad \forall u,v\in V,
    \label{eq:operatorAN}
\end{equation}
and
\begin{equation}
    (K^N(k)u,v)_{H^1(\Omega)}
    =-((k^2n(x)+1)u,v) 
    -\frac{\pi}{2}\sum_{n=1}^{N} kR\frac{H_{n-1}^{(1)}(kR)}{H_n^{(1)}(kR)} u_n\overline{v_n} \quad \forall u,v\in V.
    \label{eq:operatorKN}
\end{equation}

\begin{lemma}\label{lem:BN}
The operator functions $B^N(\cdot): \Lambda\to\calL(V)$ are holomorphic and Fredholm with index zero. Furthermore, for each compact $\Lambda_1\subset \Lambda$, $\| B^N(k)\|\le C$, $\forall N$, $\forall k\in \Lambda_1$.
\end{lemma}

\begin{proof}
As in the proof of Theorem \ref{thm:B}, one can show that $B^N(\cdot)$ is holomorphic and Fredholm with index zero on $\Lambda$.
Furthermore, by \eqref{eq:TNbound}, for each compact $\Lambda_1\subset \Lambda$,
\begin{equation*}
    \|T^N(k)\varphi\|_{H_D^{-1/2}(\Gamma^+_R)}\le
    C\|\varphi\|_{H_D^{1/2}(\Gamma^+_R)}, \quad k\in\Lambda_1,
\end{equation*}
where $C$ depends only on $\Lambda_1$. Therefore, by \eqref{eq:operatorBN},
\begin{equation*}
    \| B^N(k)\| 
    \le C,\quad\forall N, \ \forall k\in\Lambda_1.
\end{equation*}
\end{proof}

Let $\Lambda_0$ be a compact subset of $\Lambda$ so that $\partial \Lambda_0\subset \rho(B)$ and $\Lambda_0\cap \sigma(B) = \{k_0\}$.
We apply Theorem \ref{thm:Karma1} to obtain a spectral approximation
of $B^N(\cdot)$ to $B(\cdot)$, where we take $X=Y=X_n=Y_n=V$,
and the operators $p_n$ and $q_n$ are identities.

\begin{theorem}\label{thm:converge-kN}
Assume that $\gamma_0 G(B,k_0) |_{\Gamma_R^+} \subset H_D^{1/2+s}(\Gamma_R^+)$ with $s>0$. For sufficiently large $N$, it holds that $\nu(\Lambda_0,B) = \nu(\Lambda_0,B^N)$ and, for each $k^N\in \sigma(B^N)\cap\Lambda_0$, 
\begin{equation}
    |k_0 - k^N|\le \frac{C}{N^{s/\kappa}}.
\end{equation}
Here $\kappa$ is the maximum length of all the Jordan chains of $B(k_0)$.
\end{theorem}

\begin{proof}
Analogous to the proof of Theorem \ref{thm:B}, we can show that $A^N$ is bounded and coercive with $\| A^N u\|_{H^1(\Omega)}\ge \|u\|_{H^1(\Omega)}$ for all $N$. Hence $\| (A^N)^{-1}\|\le C$ for all $N$. By the definitions \eqref{eq:operatorA} and \eqref{eq:operatorAN},
\begin{align*}
    |((A-A^N)u,v)_{H^1(\Omega)}| =& \left|\frac{\pi}{2}\sum_{n=N+1}^{\infty}nu_n\overline{v_n}\right| \\
    \le{}& C\gamma_{1/2}(N,u)\|u\|_{H_D^{1/2}(\Gamma^+_R)}\|v\|_{H_D^{1/2}(\Gamma^+_R)} \\
    \le{}& C\gamma_{1/2}(N,u)\|u\|_{H^{1}(\Omega)}\|v\|_{H^{1}(\Omega)}\quad \forall u, v\in V.
\end{align*}
Therefore $A^N$ converges stably to $A$.
By the definitions \eqref{eq:operatorK}, \eqref{eq:operatorKN}
and the estimate \eqref{eq:hdh2},
\begin{align*}
    |((K(k)-K^N(k))u,v)_{H^1(\Omega)}| 
    =& \left|\frac{\pi}{2}\sum_{n=N+1}^{\infty}kR\frac{H^{(1)}_{n-1}(kR)}{H^{(1)}_n(kR)}u_n\overline{v_n}\right| \\
    \le{}& \frac{C}{N^2}\|u\|_{H_D^{1/2}(\Gamma^+_R)}\|v\|_{H_D^{1/2}(\Gamma^+_R)} \\
    \le{}& \frac{C}{N^2}\|u\|_{H^{1}(\Omega)}\|u\|_{H^{1}(\Omega)} \quad \forall u, v\in V.
\end{align*}
Therefore $K^N(k)$ converges uniformly to $K(k)$. Together with the stably convergence of $A^N$ to $A$, Lemma \ref{lem:reg-conv} implies the regular convergence of $B^N(k)$ to $B(k)$. Then by Lemma \ref{lem:BN} and Theorem \ref{thm:Karma1}, it holds that $\nu(\Lambda_0,B) = \nu(\Lambda_0,B^N)$ and
\begin{equation}
    |k_0 - k^N| 
    \le C\sup_{\substack{u\in G(B,k_0)\\ \| u\|_{H^1(\Omega)}=1, k\in\partial\Lambda_0}}
    \|(B(k)-B^N(k))u\|_{H^1(\Omega)}^{1/\kappa}.
\end{equation}

By the definitions \eqref{eq:operatorB} and \eqref{eq:operatorBN}, it holds that
\begin{align}
    & |((B(k)-B^N(k))u,v)_{H^1(\Omega)}| \notag \\
    ={}&  |\langle (T(k)-T^N(k))u,v\rangle_{\Gamma_R^+}| \notag\\
    \le{}& C\|(T(k)-T^N(k))u\|_{H_D^{-1/2}(\Gamma^+_R)}\| v\|_{H^1(\Omega)} \quad \forall u, v\in V.
\end{align}
Lemma \ref{lem:Truncation} implies that
\begin{equation}
    \|(T(k)-T^N(k))u\|_{H_D^{-1/2}(\Gamma^+_R)}
    \le \frac{C}{N^s}\| u\|_{H_D^{1/2+s}(\Gamma_R^+)} \quad \forall u\in H_D^{1/2+s}(\Gamma_R^+).
\end{equation}
Using the assumption $\gamma_0 G(B,k_0) |_{\Gamma_R^+} \subset H_D^{1/2+s}(\Gamma_R^+)$ and the finite dimension of $G(B,k_0)$, we have that
\begin{equation*}
    |k_0 - k^N|
    \le C\sup_{\substack{u\in G(B,k_0)\\ \| u\|_{H^1(\Omega)}=1, k\in\partial\Lambda_0}}
    \|(T(k)-T^N(k))u\|_{H_D^{-1/2}(\Gamma^+_R)}^{1/\kappa}
    \le \frac{C}{N^{s/\kappa}}.
\end{equation*}
The proof is complete.
\end{proof}

\section{Finite element method}\label{FEM}

Now we present a finite element method for $B^N(\cdot)$ \cite{SunZhou2016}. Let $V_n := V_{h_n} \subset V$ be some conforming finite element space with mesh size $h_n$, where $h_n \to 0$ as $n \to \infty$. Such space can be constructed using curved elements as follows. First, generate a triangular mesh $\tilde{\mathcal{T}}_{h_n}$ for $\Omega$ with mesh size $h_n$. To exactly cover $\Omega$, construct a curvilinear partition $\mathcal{T}_{h_n}$ by deforming each element $\tilde{K} \in \tilde{\mathcal{T}}_{h_n}$ into a curvilinear element $K \in \mathcal{T}_{h_n}$ via a mapping $\tilde{F}_{\tilde{K}}:\tilde{K}\to K$.
The conforming finite element space $V_n$ is defined by
\begin{equation}\label{Vn}
    V_n = \{u_h\in V\,:\, u_h|_K\circ \tilde{F}_{\tilde{K}} = p_h, \, \text{for some\,} \, p_h\in P_1(\tilde{K}),\; \tilde{K}\in \tilde{\mathcal{T}}_{h_n}\},
\end{equation}
where $P_1$ is the linear polynomial space.
Standard approximation property holds for this space $V_n$.
See \cite{Dubios} for the details of the curved element space.

Define the operator $B^N_n(k) : V_n \to V_n$ such that
\begin{equation}
    \label{eq:operatorBNn}
    (B^N_n(k)u,v)_{H^1(\Omega)} = (B^N(k)u,v)_{H^1(\Omega)}\quad \forall u,v \in V_n.
\end{equation}
The finite element approximation to \eqref{eq:eigN} is defined by
\begin{equation}
    \label{eq:eigNn}
    B^N_n(k)u^N_n=0 \quad \text{in } V_n.
\end{equation} 

\begin{lemma}\label{lem:BNn}
The resolvent set $\rho(B^N)$ is non-empty for $N$ large enough. The operator functions $B^N_n(\cdot): \Lambda\to\calL(V_n)$ are holomorphic and Fredholm with index zero. Furthermore, for each compact $\Lambda_1\subset \Lambda$, $\| B^N_n(k)\|\le C$, $\forall N$, $\forall n$, $\forall k\in \Lambda_1$.
\end{lemma}

\begin{proof}
By Theorem \ref{thm:converge-kN}, $\nu(\Lambda_0,B) = \nu(\Lambda_0,B^N)$, we know that $\rho(B^N)$ is non-empty.
The rest statement can be proved analogously as Lemma \ref{lem:BN}.
\end{proof}
Take $X=Y=V$ and $X_n=Y_n=V_n$.
Let $q_n : V\to V_n$ be the $H^1$-projection
\begin{equation}
    (q_n u,v)_{H^1(\Omega)} = (u,v)_{H^1(\Omega)}\quad \forall u\in V, \;\forall v\in V_n,
\end{equation}
and let $r_n : V_n\to V$ be the embedding operator.
Then it hold that
\begin{equation}\label{eq:BNn}
    B^N_n(k) = q_n B^N(k) r_n.
\end{equation}
Similar to the continuous case, we split $B^N_n(k)$ by
\begin{equation}
    B^N_n(k) = A^N_n + K^N_n(k)
\end{equation}
with $A^N_n : V_n\to V_n$ and $K^N_n(k) : V_n\to V_n$ defined respectively by
\begin{equation}\label{eq:ANn}
    A^N_n = q_n A^N r_n\quad \text{and} \quad K^N_n(k) = q_n K^N_n(k) r_n.
\end{equation}
By the approximation property of $V_n$, it holds that
\begin{equation}\label{eq:proj}
    \|u-q_n u\|_{H^1(\Omega)}\rightarrow 0,\quad n\rightarrow \infty,\ \ \forall u\in V.
\end{equation}
Let $p_n = q_n$, then the operators $p_n$ and $q_n$ satisfy the conditions for the spectral approximation theory of Appendix \ref{sec2.3}.

\begin{theorem}
Let $\Lambda_0$ be a compact subset of $\Lambda$ so that $\partial \Lambda_0\subset \rho(B)$ and $\Lambda_0\cap \sigma(B) = \{k_0\}$. If $G(B^N,k^N)\subset H^{1+s_N}(\Omega)\cap V$ with $s_N>0$, then, for sufficiently large $N$ and $n$, it holds that $\nu(\Lambda_0,B^N) = \nu(\Lambda_0,B^N_n)$ and, for each $k^N_n\in \sigma(B^N_n)\cap \Lambda_0$,
\begin{equation}
    |k^N - k^N_n|\le C h_n^{2s_N/\kappa_N}.
\end{equation}
 Here $\kappa_N$ is the maximum length of all the Jordan chains of $B^N(\cdot)$ at $k^N$.
\end{theorem}

\begin{proof}
Due to the proof of Theorem \ref{thm:B}, $\|(A^N_n)^{-1}\|\le C$ for all $N$ and $n$. By \eqref{eq:ANn}, we have that
\begin{align*}
    \| q_n A^N u - A^N_n p_n u\|_{H^1(\Omega)} ={}& \| q_n A^N (u-p_n u)\|_{H^1(\Omega)} \\
    \le{}&\|q_n\| \|A^N\| \|u - p_n u\|_{H^1(\Omega)}\to 0.
\end{align*}
Consequently, $A^N_n$ converges stably to $A^N$. For $0<\epsilon<1/2$ and for all $u\in V$, we have the boundedness of $K^N(k) : H^{1/2+\epsilon}(\Omega)\cap V\to V$:
\begin{equation}\label{eq:KNbound}
    \|K^N(k)u\|_{H^1(\Omega)}\leqslant C\|u\|_{L^2(\Omega)} + C\|u\|_{H_D^{-3/2}(\Gamma_R^+)} \leqslant C\|u\|_{H^{1/2+\epsilon}(\Omega)}.
\end{equation}
By the equalities \eqref{eq:ANn} and \eqref{eq:KNbound}, we have that
\begin{align*}
    &\|q_n K^N(k) u - K^N_n(k) p_n u\|_{H^1(\Omega)} \\ ={}& \| q_n K^N(k) (u-p_n u)\|_{H^1(\Omega)} \\
    \le{}&\|q_n\| \|K^N(k)\|_{H^{1/2+\epsilon}(\Omega)\to H^1(\Omega)} \|u - p_n u\|_{H^{1/2+\epsilon}(\Omega)} \\
    \le{}& C h_n^{1/2-\epsilon}\| u\|_{H^1(\Omega)}.
\end{align*} 
Thus $K^N_n(k)$ converges uniformly to $K^N(k)$. By Lemma \ref{lem:reg-conv}, $B^N_n(k)$ converges regularly to $B^N(k)$. Hence Lemma \ref{lem:BNn} and Theorem \ref{thm:Karma1} imply that $\nu(\Lambda_0,B^N) = \nu(\Lambda_0,B^N_n)$.

Next, the operator $r_n=q_n^*$ satisfy the assumptions of Theorem \ref{thm:Karma2}, moreover \eqref{eq:BNn} holds. Similar to \eqref{eq:Tuv}, we have that
\begin{equation}
    \langle T^N(k)u,v\rangle_{\Gamma_R^+} = \frac{\pi R}{2}\sum_{n=1}^N \tau_n(k)u_n\overline{v_n} = \langle T^N(k)\overline{v},\overline{u}\rangle_{\Gamma_R^+}.
\end{equation}
By the definition \eqref{eq:operatorBN}, it holds that
\begin{equation}
    (B^N(k)u,v)_{H^1(\Omega)} = (B^N(k)\overline{v},\overline{u})_{H^1(\Omega)}.
\end{equation}
Therefore, 
\[
(B^N)^*(k)v := [B^N(\overline{k})]^*v = \overline{B^N(\overline{k})\overline{v}}
\] 
and, for each $j$, 
\[
((B^N)^*)^{(j)}(k)v = \overline{(B^N)^{(j)}(\overline{k})\overline{v}}.
\]
Consequently, $G((B^N)^*,\overline{k^N}) = \overline{G(B^N,k^N)}$ and
\begin{align}
    d^*_n 
    &:= \max_{\substack{u\in G((B^N)^*,\overline{k^N})\\ \|u\|_{H^1(\Omega)} = 1}} \textnormal{dist}(u,q_n^* V_n)
    = \max_{\substack{w\in G(B^N,k^N)\\ \|w\|_{H^1(\Omega)} = 1}} \textnormal{dist}(\overline{w},r_n V_n) \notag\\
    &\;= \max_{\substack{w\in G(B^N,k^N)\\ \|w\|_{H^1(\Omega)} = 1}} \textnormal{dist}(w,r_n V_n)
    := d_n.
\end{align}
Thus
\begin{equation}
    d^*_n = d_n 
    = \max_{\substack{u\in G(B^N,k^N) \\ \|u\|_{H^1(\Omega)}=1}} 
    \min_{v\in V_n} \| u-v\|_{H^1(\Omega)}
    \le C h_n^{s_N},
\end{equation}
where we have used the assumption $G(B^N,k^N)\subset H^{1+s_N}(\Omega)$ and the finite dimension of $G(B^N,k^N)$. Application of Theorem \ref{thm:Karma2} completes the proof.
\end{proof}

\section{Numerical examples}\label{NE}

In this section, we present numerical examples to demonstrate the effectiveness of the proposed method. The linear finite elements are used on a sequence of uniformly refined meshes $\{\calT_{h_l}\}_l$, which may have curved boundaries along $\partial\Omega$. The truncation order is $N = 10$ for all examples since increasing $N$ beyond 10 leads to only minor changes in the computed resonances.

To solve the nonlinear eigenvalue problem \eqref{eq:eigNn}, we apply the spectral indicator method (SIM) for the search region
\begin{equation*}
    \Theta = \{ a+b\rmi \in \bbC : a\in[0,L],\ b\in[-L,0] \},
\end{equation*}
where $L=6$ for the first three examples and $L=8$ for the last example. SIM is a contour-integral-based method originally developed for non-Hermitian matrix eigenvalue problems~\cite{Huang2016JCP,Huang2018NLAA}. Highly scalable parallel implementations of SIM have also been successfully applied to nonlinear eigenvalue problems~\cite{Gong2022MC,MaSun2023,XiLinSun2024CMA}. Since SIM requires no prior estimates of the resonance locations, it is well suited to solving complex eigenvalue problems.

To evaluate the numerical results, we define the relative error and the convergence order, respectively, as
\begin{equation*}
    E_j= \frac{\left|k^j - k^{j+1}\right|}{\left|k^{j+1}\right|},
    \qquad j = 1,2,3,4,
\end{equation*}
and
\begin{equation*}
    \log_2\left(\frac{E_j}{E_{j+1}}\right),
    \qquad j = 1,2,3,
\end{equation*}
where \(k^j\) denotes the resonance computed on the mesh obtained after the \(j\)-th refinement.

\subsection{Overfilled circular cavity}

\begin{figure}[htbp!]
    \centering
    \begin{subfigure}{0.48\linewidth}
        \centering
        \includegraphics[width=\linewidth]{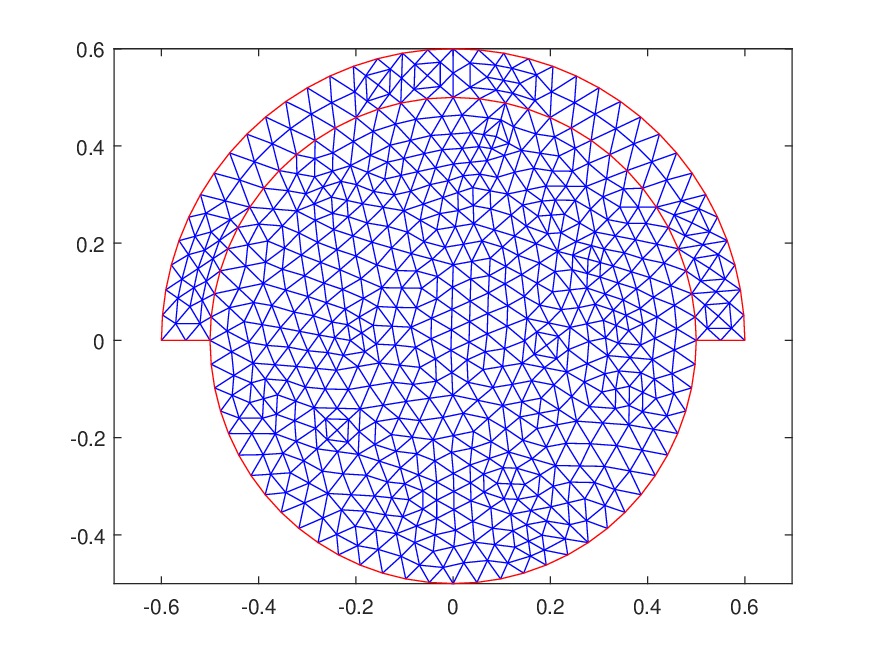}
        \caption{}
        \label{fig:circle:a}
    \end{subfigure}
    \hfill
    \begin{subfigure}{0.48\linewidth}
        \centering
        \includegraphics[width=\linewidth]{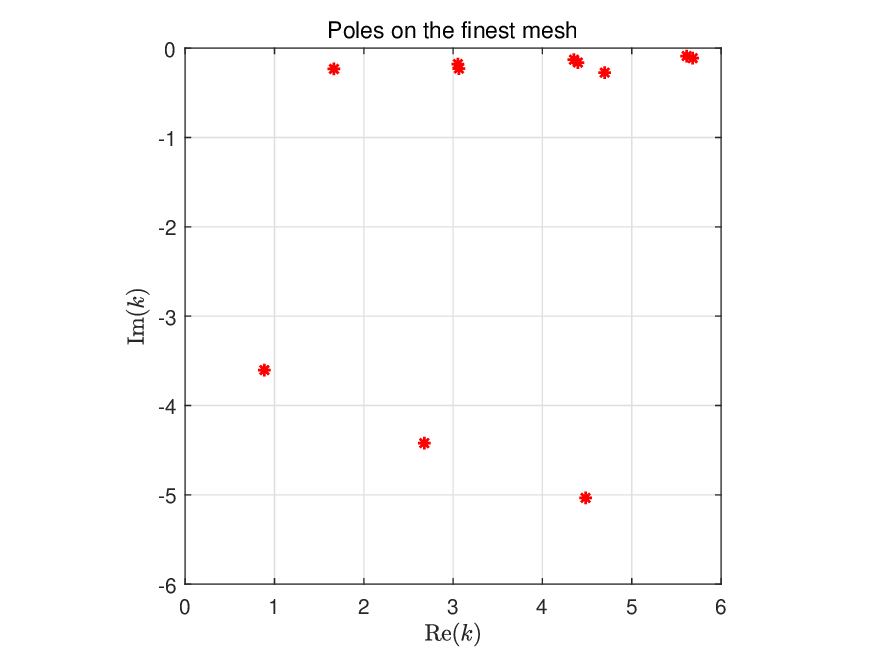}
        \caption{}
        \label{fig:circle:b}
    \end{subfigure}
    \caption{Overfilled circular cavity (a) and the computed resonances (b).}
    \label{fig:circle}
\end{figure}

We consider a circular cavity with radius $r=0.5$. Half of the cavity is under the ground plane. The cavity is filled with a homogeneous medium whose refractive index is $n(x)=4$. The computational domain is terminated by a half circle with radius $R=0.6$, on which the DtN map is imposed. Figure \ref{fig:circle:a} shows the geometry of the cavity and the initial mesh of the computational domain.

Figure \ref{fig:circle:b} shows the computed resonances in $\Theta$. Table \ref{tab:circle} lists four resonances and their convergence orders. The first two resonances lie close to the real axis, whereas the other two lie away from it. The convergence orders are approximately 1.5 for the first two resonances and 2 for the other two.

\begin{table}[htbp!]
    \caption{Convergence orders for overfilled circular cavity.}
    \label{tab:circle}
    \centering
    \resizebox{\linewidth}{!}{
        \begin{tabular}{lcccccccc}
            \hline
            &$k_{1,h}$&Ord&$k_{2,h}$&Ord&$k_{3,h}$&Ord&$k_{4,h}$&Ord\\
            \hline
            DOF = 181&1.6826 - 0.2366i&&3.1022 - 0.1874i&&0.8817 - 3.5701i&&2.6321 - 4.3132i&\\
            \hline
            DOF = 684&1.6714 - 0.2338i&&3.0676 - 0.1800i&&0.8861 - 3.5944i&&2.6654 - 4.3929i&\\
            \hline
            DOF = 2665&1.6676 - 0.2330i&1.56&3.0568 - 0.1780i&1.68&0.8872 - 3.6016i&1.78&2.6753 - 4.4154i&1.82\\
            \hline
            DOF = 10513&1.6662 - 0.2328i&1.49&3.0532 - 0.1774i&1.60&0.8874 - 3.6037i&1.81&2.6778 - 4.4214i&1.93\\
            \hline
            DOF = 41761&1.6657 - 0.2327i&1.44&3.0520 - 0.1772i&1.53&0.8874 - 3.6043i&1.76&2.6783 - 4.4230i&1.95\\
            \hline
        \end{tabular}
    }
\end{table}

\subsection{Overfilled cavity with an enlarged semicircle}

\begin{figure}[htbp]
    \centering
    \begin{subfigure}{0.48\linewidth}
        \centering
        \includegraphics[width=\linewidth]{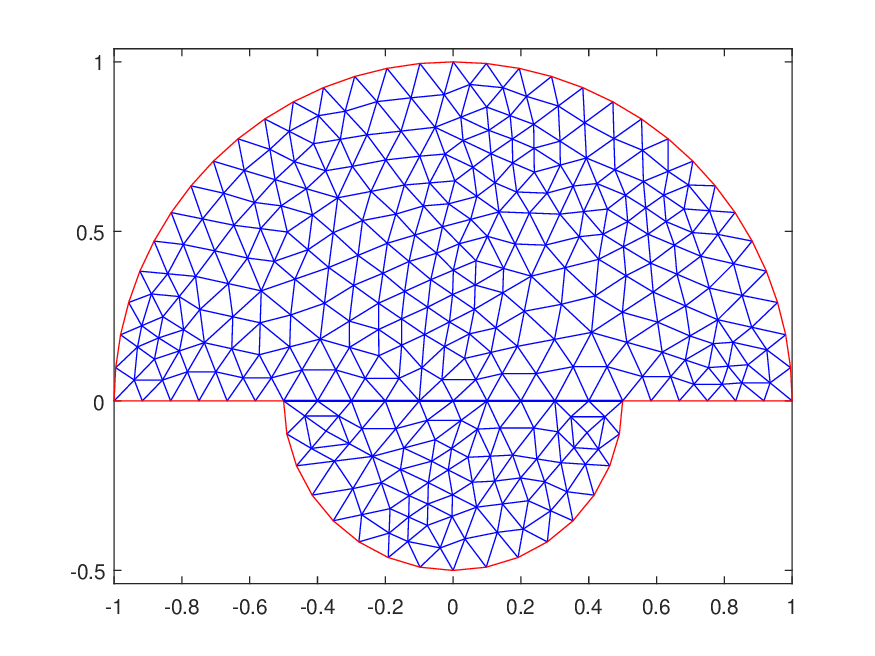}
        \caption{}
        \label{fig:circle-large:a}
    \end{subfigure}
    \begin{subfigure}{0.48\linewidth}
        \centering
        \includegraphics[width=\linewidth]{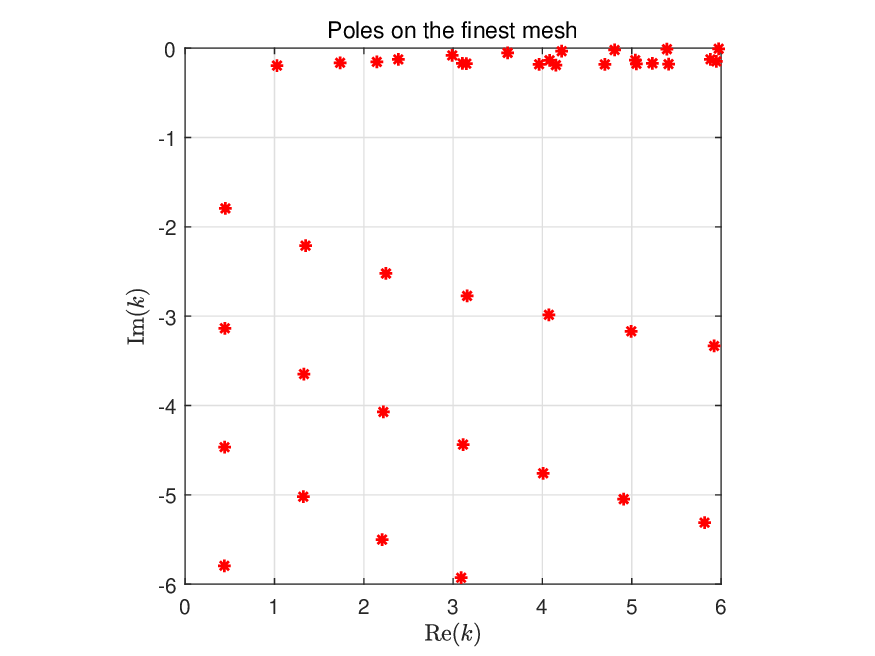}
        \caption{}
        \label{fig:circle-large:b}
    \end{subfigure}
    \caption{Overfilled cavity (a) and the computed resonances (b).}
    \label{fig:circle-large}
\end{figure}

In this example, the cavity consists a half disk with $r=0.5$ under the ground plane and a half disk with $r=1.0$ above the ground plane (see Fig. \ref{fig:circle-large:a}). It is again filled with a homogeneous medium with constant refractive index $n(x)=4$. The DtN boundary condition is exactly imposed on the top semicircle.

Fig.~\ref{fig:circle-large:b} shows the computed resonances in $\Theta$. Table~\ref{tab:circle-large} lists four resonances together with their convergence orders. The first two resonances exhibit a convergence order of approximately 1.5, whereas the last two exhibit an order of approximately 2.

\begin{table}[htbp!]
    \caption{Convergence orders of the resonances (Example 2).}
    \label{tab:circle-large}
    \begin{center}
        \resizebox{\linewidth}{!}{
            \begin{tabular}{lcccccccc}
                \hline
                &$k_{1,h}$&Ord&$k_{2,h}$&Ord&$k_{3,h}$&Ord&$k_{4,h}$&Ord\\
                \hline
                DOF = 93&1.0425 - 0.2022i&&1.7598 - 0.1757i&&0.4483 - 1.7836i&&1.3395 - 2.1897i&\\
                \hline
                DOF = 339&1.0331 - 0.1969i&&1.7418 - 0.1668i&&0.4499 - 1.7909i&&1.3453 - 2.2056i&\\
                \hline
                DOF = 1293&1.0298 - 0.1952i&1.55&1.7365 - 0.1643i&1.77&0.4504 - 1.7929i&1.90&1.3471 - 2.2099i&1.88\\
                \hline
                DOF = 5049&1.0287 - 0.1946i&1.49&1.7349 - 0.1636i&1.71&0.4505 - 1.7934i&1.96&1.3475 - 2.2109i&1.97\\
                \hline
                DOF = 19953&1.0282 - 0.1944i&1.44&1.7344 - 0.1633i&1.64&0.4505 - 1.7936i&1.98&1.3477 - 2.2112i&2.00\\
                \hline
            \end{tabular}
        }
    \end{center}
\end{table}

\subsection{Overfilled cavity with a narrow channel}

\begin{figure}[htbp]
    \centering
    \begin{subfigure}{0.48\linewidth}
        \centering
        \includegraphics[width=\linewidth]{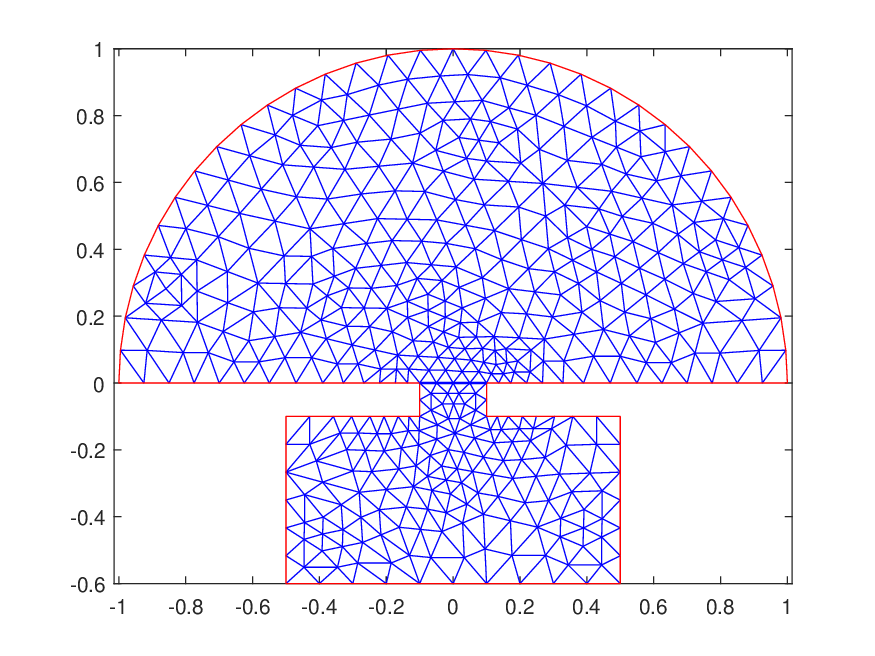}
        \caption{}
        \label{fig:channel:a}
    \end{subfigure}
    \begin{subfigure}{0.48\linewidth}
        \centering
        \includegraphics[width=\linewidth]{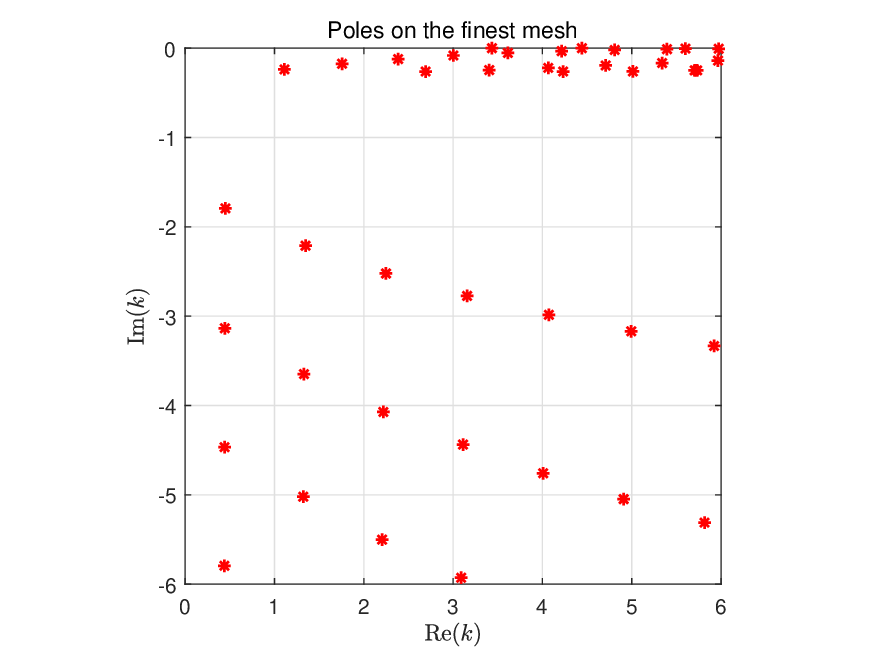}
        \caption{}
        \label{fig:channel:b}
    \end{subfigure}
    \caption{Cavity with a narrow channel and the computed resonances.}
    \label{fig:cave}
\end{figure}

The top part of the cavity is a half disk with $r=1.0$ above the ground plane. The lower part is 
\begin{equation*}
    \Omega_d = ([-0.1, 0.1]\times[-0.1, 0])\cup([-0.5, 0.5]\times[-0.6, -0.1]),
\end{equation*}
where a narrow channel exists (Fig. \ref{fig:channel:a}). The entire cavity is filled with a homogeneous medium with refractive index $n(x)=4$.

Fig.~\ref{fig:channel:b} shows the computed resonances in  $\Theta$.
Table~\ref{tab:channel} lists four resonances together with their convergence orders. All four resonances exhibit approximately second-order convergence.

\begin{table}[htbp]
    \caption{Convergence orders of the resonances (Example 3).}
    \label{tab:channel}
    \begin{center}
        \resizebox{\linewidth}{!}{
            \begin{tabular}{ccccccccc}
                \hline
                &$k_{1,h}$&Ord&$k_{2,h}$&Ord&$k_{3,h}$&Ord&$k_{4,h}$&Ord\\
                \hline
                DOF = 134&1.1191 - 0.2438i&&1.7720 - 0.1829i&&0.4483 - 1.7835i&&1.3395 - 2.1881i&\\
                \hline
                DOF = 482&1.1138 - 0.2393i&&1.7603 - 0.1765i&&0.4502 - 1.7908i&&1.3453 - 2.2053i&\\
                \hline
                DOF = 1823&1.1123 - 0.2381i&1.84&1.7573 - 0.1749i&1.95&0.4507 - 1.7927i&1.90&1.3471 - 2.2098i&1.88\\
                \hline
                DOF = 7085&1.1119 - 0.2378i&1.80&1.7565 - 0.1745i&1.98&0.4508 - 1.7932i&1.97&1.3476 - 2.2110i&1.96\\
                \hline
                DOF = 27929&1.1117 - 0.2377i&1.74&1.7563 - 0.1744i&1.99&0.4509 - 1.7933i&1.99&1.3477 - 2.2113i&1.99\\
                \hline
            \end{tabular}
        }
    \end{center}
\end{table}

\subsection{Example 4}

\begin{figure}[htbp]
    \centering
    \begin{subfigure}{0.48\linewidth}
        \centering
        \includegraphics[width=\linewidth]{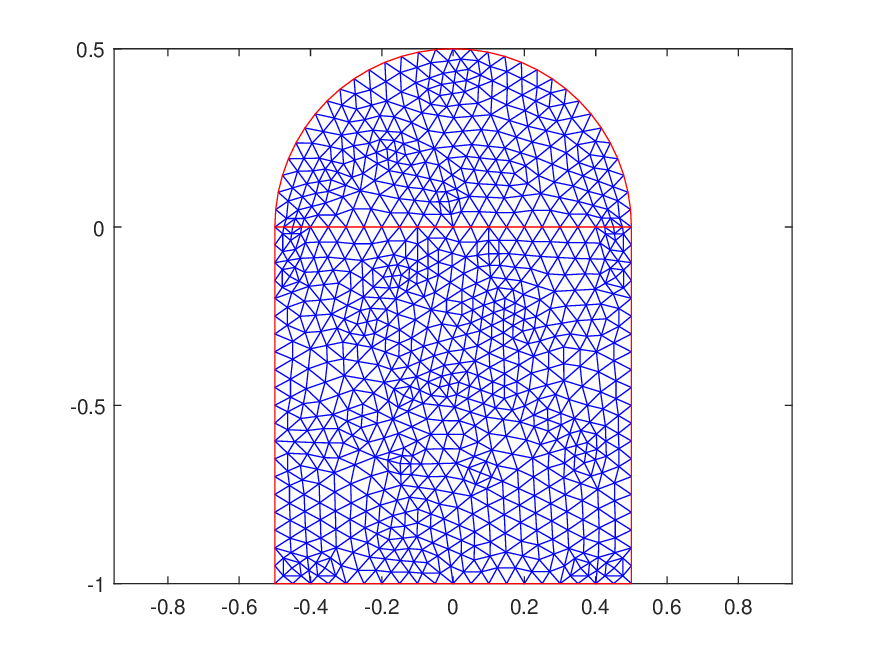}
        \caption{}
        \label{fig:square:a}
    \end{subfigure}
    \hfill
    \begin{subfigure}{0.48\linewidth}
        \centering
        \includegraphics[width=\linewidth]{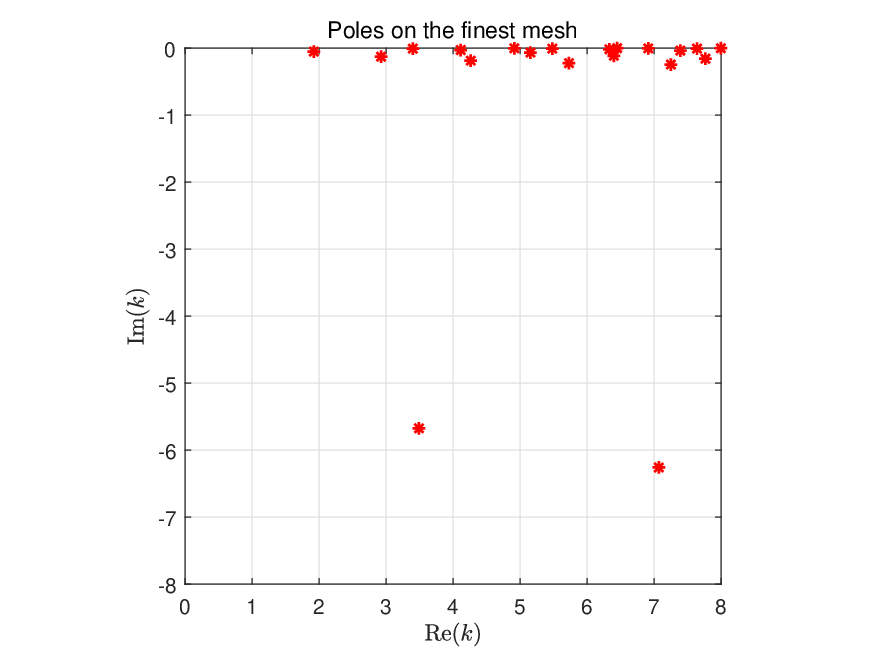}
        \caption{}
        \label{fig:square:b}
    \end{subfigure}
    \caption{Unit-square cavity and the computed resonances.}
    \label{fig:square}
\end{figure}

The cavity is the unit square under the ground plane (see Fig.~\ref{fig:square:a}). It is filled with a homogeneous medium with refractive index $n(x)=4$. The computational domain is truncated by an upper semi-circle of radius $r=0.5$, on which the DtN map is imposed.

Fig.~\ref{fig:square:b} shows the computed resonances in $\Theta$. Table~\ref{tab:square} lists four  resonances together with their convergence orders, which are close to 1.5 for all the resonances.

\begin{table}[htbp]
    \caption{Convergence orders of selected resonances for unit-square cavity.}
    \label{tab:square}
    \centering
    \resizebox{\linewidth}{!}{
        \begin{tabular}{ccccccccc}
            \hline
            &$k_{1,h}$&Ord&$k_{2,h}$&Ord&$k_{3,h}$&Ord&$k_{4,h}$&Ord\\
            \hline
            DOF = 244&1.9318 - 0.0494i&&2.9541 - 0.1282i&&3.4243 - 4.8711i&&7.4220 - 6.6018i&\\
            \hline
            DOF = 927&1.9244 - 0.0499i&&2.9325 - 0.1264i&&3.5626 - 5.3343i&&6.6905 - 5.6240i&\\
            \hline
            DOF = 3613&1.9220 - 0.0502i&1.60&2.9263 - 0.1264i&1.82&3.5372 - 5.5977i&0.91&6.9498 - 6.0360i&1.40\\
            \hline
            DOF = 14265&1.9210 - 0.0505i&1.38&2.9245 - 0.1266i&1.69&3.4904 - 5.6744i&1.57&7.0700 - 6.2588i&0.98\\
            \hline
            DOF = 56689&1.9207 - 0.0506i&1.24&2.9238 - 0.1267i&1.52&3.4694 - 5.6870i&1.88&7.1022 - 6.3393i&1.56\\
            \hline
        \end{tabular}
    }
\end{table}

\appendix
\section{Spectral Approximation Theory}\label{sec2.3}
We present the main results of the spectral approximation of eigenvalues of holomorphic Fredholm operator functions \cite{Vnk76, Karma1, Karma2}. Let $X$, $Y$, $X_n$, $Y_n$ be Banach spaces and let the operators $p_n\in\mathcal{L}(X,X_n)$ and $q_n\in \mathcal{L}(Y,Y_n)$ satisfy $\| p_n x\|_{X_n} \to \| x\|_{X}$ and $\| q_n y\|_{Y_n} \to \| y\|_{Y}$ as $n\to \infty$ for all $x\in X$ and $y\in Y$. 

Let $\{x_n\}$ be such that $x_n \in X_n$, we define compactness of $\{x_n\}$ as follows. 
If for every subsequence $\{x'_n\}$ of $\{x_n\}$ there exists a subsequence $\{x''_n\}$ of $\{x'_n\}$ and $x\in X$ such that 
$\| x''_n - p_n x\| \to 0$ as $n\to \infty$, then $\{x_n\}$ is called (discrete) compact. The compactness of $\{y_n\}$ with $y_n\in Y_n$ is defined analogously.

Let $F\in \mathcal{L}(X,Y)$ and $F_n\in \mathcal{L}(X_n,Y_n)$.
Definitions of convergence of $F_n$ to $F$ are given below.
\begin{itemize}
\item
$F_n$ is said to converge to $F$, if $\| F_n p_n x - q_n F x\| \to 0,\; n\to\infty,\; \forall x\in X$.
\item
$F_n$ is said to converge uniformly to $F$, if $\| F_n p_n - q_n F\|\to 0$, $n\to \infty$.
\item
$F_n$ is said to converge stably to $F$, if $F_n$ converges to $F$ and there exists $n_0$ such that $F_n$ is invertible with $\| F_n^{-1}\|\le C$, $\forall n\ge n_0$.
\item
$F_n$ is said to converge compactly to $F$, if $F_n$ converges to $F$ and $\{F_n x_n\}$ is compact for any $\| x_n\| \le 1$, $\forall n$.
\item
$F_n$ is said to converge regularly to $F$, if $F_n$ converges to $F$ and $\{x_n\}$ is compact for any $\| x_n\|\le 1$, $\forall n$ such that
$\{F_n x_n\}$ is compact. 
\end{itemize}

Let $E, H\in \mathcal{L}(X,Y)$ and $E_n, H_n \in \mathcal{L}(X_n,Y_n)$, then the following Lemma can be used to prove regular convergence (\cite[Lemma 3.4]{XiGongSun2024}).

\begin{lemma}\label{lem:reg-conv}
If $E$ is invertible, $H$ is compact, $E_n$ converges stably to $E$ and $H_n$ converges uniformly to $H$, then
$E_n+H_n$ converges regularly to $E+H$.
\end{lemma}

Let $\Lambda_0$ be a compact subset 
of $\Lambda$ so that $\partial\Lambda_0\subset
\rho(F)$ and $\Lambda_0\cap \sigma(F) = \{\lambda_0\}$. Denote by $\nu(\Lambda_0,F)$ (or $\nu(\Lambda_0,F_n)$) the sum of algebraic multiplicities of all the eigenvalues
of $F(\cdot)$ (or $F_n(\cdot)$) in $\Lambda_0$.
Then by Theorems 2 and 3 of \cite{Karma1}
and Theorem~2 of \cite{Karma2},
the following result holds.

\begin{theorem}\label{thm:Karma1}
Let $F(\cdot) : \Lambda\to \mathcal{L}(X,Y)$
be a holomorphic Fredholm operator function with index zero and $\rho(F)\neq \emptyset$.
Let $F_n(\cdot) : \Lambda\to \mathcal{L}(X_n,Y_n), n\in\bbN$
be a sequence of holomorphic Fredholm operator functions with index zero.
Assume that $\{F_n(\cdot)\}_{n\in\bbN}$ is equibounded on any compact $\Lambda_1\subset \Lambda$, i.e.,
\begin{equation*}
\| F_n(\lambda) \| \le C,\quad \forall \lambda\in \Lambda_1,\; n\in\bbN,
\end{equation*}
and $F_n(\lambda)$ converges regularly to $F(\lambda)$ for every $\lambda\in \Lambda$.
Then for every compact $\Lambda_0\subset \Lambda$ with $\partial\Lambda_0\subset\rho(F)$ and for sufficiently large $n$,
it holds that $\nu(\Lambda_0,F_n) = \nu(\Lambda_0,F)$.
Moreover, for each $\lambda_n\in \sigma(F_n)\cap\Lambda_0$,
\begin{equation*}
|\lambda_0 - \lambda_n|\le C\sup_{\substack {x\in G(F,\lambda_0)\\
\| x\| = 1,\; \lambda\in \partial\Lambda_0}}\| q_n F(\lambda)x - F_n(\lambda)p_n x\|^{1/\kappa}.
\end{equation*}
Here $\kappa$ is the maximum length of all the Jordan chains of $F(\cdot)$ at $\lambda_0$.
\end{theorem}

\begin{theorem}\label{thm:Karma2} (\cite[Theorem 3]{Karma2})
Let the assumptions of Theorem \ref{thm:Karma1} hold.
Let $r_n\in \mathcal{L}(X_n,X)$ and $q_n\in \mathcal{L}(Y,Y_n)$ be such that
\begin{align*}
    &\| r_n\| \le C,\; \forall n, \\
    &\| r_np_n x - p_n x\|_{X_n} \to 0,\; n\to \infty, \; \forall x\in G(F,\lambda_0).
    \end{align*}
If $F_n(k) = q_n F(k)r_n$ for all $k\in \Lambda$, then for sufficiently large $n$ and each $\lambda_n \in \sigma(F_n)\cap \Lambda_0$, it holds that
\begin{equation*}
    |\lambda_0 - \lambda_n|\le C(d_n d_n^*)^{1/\kappa},
\end{equation*}
where
\begin{align*}
d_n  = \sup_{\substack{x\in G(F,\lambda_0)\\ \| x\| = 1}} \textnormal{dist}(x,r_n X_n),\quad
d_n^* = \sup_{\substack{y\in G(F^*,\overline{\lambda_0})\\ \| y\| = 1}}\textnormal{dist}(y, q_n^*Y_n),
\end{align*}
and $\kappa$ is the maximum length of all the Jordan chains of $F(\cdot)$ at $\lambda_0$.
\end{theorem}

\section{Poles of the DtN Mapping}
The definition of $\Lambda$ excludes zeros of $H_n^{(1)}(kR)$. To compute the scattering resonances, we apply the spectral indicator method which needs to evaluate $B_N(k)$ for many quadrature points $k$'s. The DtN map is not defined at a value $k_0$ for which
\[
  H_{n}^{(1)}(k_0R)=0
\]
for at least one $n$.  These values are poles introduced by the artificial boundary and are not, in general, physical scattering resonances of the original exterior problem.

Note that the zeros of $H_{n}^{(1)}$ are simple. Let $z_0=k_0R$ is a simple zero of $H_{n}^{(1)}$ for some $n$.  Then
\[
  H_{n}^{(1)}(kR)
  =R(H_{n}^{(1)})'(k_0R)(k-k_0)+O\bigl((k-k_0)^2\bigr).
\]
Since $(H_{n}^{(1)})'(k_0R)\ne 0$, it follows that
\begin{equation}\label{eq:residue}
  \tau_n(k)
  =k\frac{(H_{n}^{(1)})'(kR)}{H_{n}^{(1)}(kR)}
  =\frac{k_0/R}{k-k_0}+O(1)
  \qquad (k\to k_0).
\end{equation}

When one use the contour-integral methods such as SIM, one can avoid using these poles as the quadrature points by checking if they are zeros of some $H_{n}^{(1)}$.

\section*{Acknowledgement}
The research of B. Gong is partially supported by National Natural Science Foundation of China No.12671468.
The research of X. Wu is partially supported by National Natural Science Foundation of China No.12671468.

%
%



\end{document}

%% file: macros.tex
\newcommand{\barB}[1]{\boldsymbol{\bar{#1}}}
\newcommand{\hatB}[1]{\boldsymbol{\hat{#1}}}
\newcommand{\tildeB}[1]{\boldsymbol{\tilde{#1}}}
\newcommand{\calB}[1]{\boldsymbol{\mathcal{#1}}}

\ExplSyntaxOn

\cs_new_protected:Npn \my_new_cmd:nnn #1#2#3
{
    \cs_if_exist:cF { #1#2 }
    {
        \cs_new_protected:cpx { #1#2 } { #3 }
    }
}

\clist_map_inline:nn {A,B,C,D,E,F,G,H,I,J,K,L,M,N,O,P,Q,R,S,T,U,V,W,X,Y,Z}
{
    \my_new_cmd:nnn {B}    {#1} { \exp_not:N \boldsymbol {#1} }
    \my_new_cmd:nnn {bB}   {#1} { \exp_not:N \barB       {#1} }
    \my_new_cmd:nnn {hB}   {#1} { \exp_not:N \hatB       {#1} }
    \my_new_cmd:nnn {tB}   {#1} { \exp_not:N \tildeB     {#1} }

    \my_new_cmd:nnn {rm}   {#1} { \exp_not:N \mathrm     {#1} }
    \my_new_cmd:nnn {V}    {#1} { \exp_not:N \mathbf     {#1} }
    \my_new_cmd:nnn {it}   {#1} { \exp_not:N \mathit     {#1} }
    \my_new_cmd:nnn {sf}   {#1} { \exp_not:N \mathsf     {#1} }
    \my_new_cmd:nnn {tt}   {#1} { \exp_not:N \mathtt     {#1} }
    \my_new_cmd:nnn {F}    {#1} { \exp_not:N \mathfrak   {#1} }

    \my_new_cmd:nnn {bb}   {#1} { \exp_not:N \mathbb     {#1} }
    \my_new_cmd:nnn {sc}   {#1} { \exp_not:N \mathscr    {#1} }
    \my_new_cmd:nnn {cal}  {#1} { \exp_not:N \mathcal    {#1} }
    \my_new_cmd:nnn {calB} {#1} { \exp_not:N \calB       {#1} }
}

\clist_map_inline:nn {a,b,c,d,e,f,g,h,i,j,k,l,m,n,o,p,q,r,s,t,u,v,w,x,y,z}
{
    \my_new_cmd:nnn {B}  {#1} { \exp_not:N \boldsymbol {#1} }
    \my_new_cmd:nnn {bB} {#1} { \exp_not:N \barB       {#1} }
    \my_new_cmd:nnn {hB} {#1} { \exp_not:N \hatB       {#1} }
    \my_new_cmd:nnn {tB} {#1} { \exp_not:N \tildeB     {#1} }

    \my_new_cmd:nnn {rm} {#1} { \exp_not:N \mathrm     {#1} }
    \my_new_cmd:nnn {V}  {#1} { \exp_not:N \mathbf     {#1} }
    \my_new_cmd:nnn {it} {#1} { \exp_not:N \mathit     {#1} }
    \my_new_cmd:nnn {sf} {#1} { \exp_not:N \mathsf     {#1} }
    \my_new_cmd:nnn {tt} {#1} { \exp_not:N \mathtt     {#1} }
    \my_new_cmd:nnn {F}  {#1} { \exp_not:N \mathfrak   {#1} }
}

\clist_map_inline:nn {Gamma,Delta,Theta,Lambda,Xi,Pi,Sigma,Upsilon,Phi,Psi,Omega}
{
    \my_new_cmd:nnn {B}  {#1} { \exp_not:N \boldsymbol { \exp_not:c {#1} } }
    \my_new_cmd:nnn {bB} {#1} { \exp_not:N \barB       { \exp_not:c {#1} } }
    \my_new_cmd:nnn {hB} {#1} { \exp_not:N \hatB       { \exp_not:c {#1} } }
    \my_new_cmd:nnn {tB} {#1} { \exp_not:N \tildeB     { \exp_not:c {#1} } }
}

\clist_map_inline:nn
{
    alpha,beta,gamma,delta,epsilon,zeta,eta,theta,iota,kappa,lambda,
    mu,nu,xi,pi,rho,sigma,tau,upsilon,phi,chi,psi,omega,
    varepsilon,vartheta,varkappa,varpi,varrho,varsigma,varphi
}
{
    \my_new_cmd:nnn {B}  {#1} { \exp_not:N \boldsymbol { \exp_not:c {#1} } }
    \my_new_cmd:nnn {bB} {#1} { \exp_not:N \barB       { \exp_not:c {#1} } }
    \my_new_cmd:nnn {hB} {#1} { \exp_not:N \hatB       { \exp_not:c {#1} } }
    \my_new_cmd:nnn {tB} {#1} { \exp_not:N \tildeB     { \exp_not:c {#1} } }
}

\ExplSyntaxOff

\DeclareMathOperator{\Curl}{Curl} 
\DeclareMathOperator{\Div}{Div} 

\newcommand{\defaultdomain}{\Omega}

\NewDocumentCommand{\LebL}{s O{\defaultdomain} m}
{
  \IfBooleanTF{#1}{\BL}{L}%
  ^{#3}(#2)%
}

\NewDocumentCommand{\SobW}{s O{\defaultdomain} m m o}
{
  \IfBooleanTF{#1}{\BW}{W}%
  \IfValueT{#5}{_{#5}}%
  ^{#3,#4}(#2)%
}

\NewDocumentCommand{\SobH}{s O{\defaultdomain} m o}
{
    \IfBooleanTF{#1}{\BH}{H}%
    \IfValueT{#4}{_{#4}}%
    ^{#3}(#2)%
}

\NewDocumentCommand{\SobHdiv}{s O{\defaultdomain} o}
{
    \IfBooleanTF{#1}{\BH}{H}%
    \IfValueT{#3}{_{#3}}%
    (\Div;\,#2)%
}

\NewDocumentCommand{\SobHcurl}{s O{\defaultdomain} o}
{
    \IfBooleanTF{#1}{\BH}{H}%
    \IfValueT{#3}{_{#3}}%
    (\Curl;\,#2)%
}

\NewDocumentCommand{\Lp}{O{\defaultdomain} m} { \LebL[#1]{#2} }
\NewDocumentCommand{\Ltwo}{O{\defaultdomain}} { \LebL[#1]{2} }
\NewDocumentCommand{\Linf}{O{\defaultdomain}} { \LebL[#1]{\infty} }

\NewDocumentCommand{\Lpv}{O{\defaultdomain} m} { \LebL*[#1]{#2} }
\NewDocumentCommand{\Ltwov}{O{\defaultdomain}} { \LebL*[#1]{2} }
\NewDocumentCommand{\Linfv}{O{\defaultdomain}} { \LebL*[#1]{\infty} }

\NewDocumentCommand{\Wmp}{O{\defaultdomain} m m} {\SobW[#1]{#2}{#3} }
\NewDocumentCommand{\bWmp}{O{\defaultdomain} m m m} { \SobW[#1]{#2}{#3}[#4] }
\NewDocumentCommand{\zbWmp}{O{\defaultdomain} m m} { \SobW[#1]{#2}{#3}[0] }

\NewDocumentCommand{\Wmpv}{O{\defaultdomain} m m} { \SobW*[#1]{#2}{#3} }
\NewDocumentCommand{\bWmpv}{O{\defaultdomain} m m m} { \SobW*[#1]{#2}{#3}[#4] }
\NewDocumentCommand{\zbWmpv}{O{\defaultdomain} m m} { \SobW*[#1]{#2}{#3}[0] }

\NewDocumentCommand{\Hm}{O{\defaultdomain} m} { \SobH[#1]{#2} }
\NewDocumentCommand{\bHm}{O{\defaultdomain} m m} { \SobH[#1]{#2}[#3] }
\NewDocumentCommand{\zbHm}{O{\defaultdomain} m} { \SobH[#1]{#2}[0] }

\NewDocumentCommand{\Hmv}{O{\defaultdomain} m} { \SobH*[#1]{#2} }
\NewDocumentCommand{\bHmv}{O{\defaultdomain} m m} { \SobH*[#1]{#2}[#3] }
\NewDocumentCommand{\zbHmv}{O{\defaultdomain} m} { \SobH*[#1]{#2}[0] }

\NewDocumentCommand{\Hone}{O{\defaultdomain}} { \SobH[#1]{1} }
\NewDocumentCommand{\bHone}{O{\defaultdomain} m} { \SobH[#1]{1}[#2] }
\NewDocumentCommand{\zbHone}{O{\defaultdomain}} { \SobH[#1]{1}[0] }

\NewDocumentCommand{\Honev}{O{\defaultdomain}} { \SobH*[#1]{1} }
\NewDocumentCommand{\bHonev}{O{\defaultdomain} m} { \SobH*[#1]{1}[#2] }
\NewDocumentCommand{\zbHonev}{O{\defaultdomain}} { \SobH*[#1]{1}[0] }

\NewDocumentCommand{\Hdiv}{O{\defaultdomain}} { \SobHdiv*[#1] }
\NewDocumentCommand{\bHdiv}{O{\defaultdomain} m} { \SobHdiv*[#1][#2] }
\NewDocumentCommand{\zbHdiv}{O{\defaultdomain}} { \SobHdiv*[#1][0] }

\NewDocumentCommand{\Hcurl}{O{\defaultdomain}} { \SobHcurl*[#1] }
\NewDocumentCommand{\bHcurl}{O{\defaultdomain} m} { \SobHcurl*[#1][#2] }
\NewDocumentCommand{\zbHcurl}{O{\defaultdomain}} { \SobHcurl*[#1][0] }

\DeclarePairedDelimiter{\norm}{\lVert}{\rVert}
\DeclarePairedDelimiter{\snorm}{\lvert}{\rvert}
\DeclarePairedDelimiter{\tnorm}{\vert\!\vert\!\vert}{\vert\!\vert\!\vert}

\DeclarePairedDelimiterX{\innerp}[2]{\lparen}{\rparen}{#1,\,#2}
\DeclarePairedDelimiterX{\dualp}[2]{\langle}{\rangle}{#1,\,#2}

\NewDocumentCommand{\NormIn}{s o m m}
{
  \IfBooleanTF{#1}
    {\norm*{#3}_{#4}}
    {
      \IfNoValueTF{#2}
        {\norm{#3}_{#4}}
        {\norm[#2]{#3}_{#4}}
    }
}

\NewDocumentCommand{\SNormIn}{s o m m}
{
  \IfBooleanTF{#1}
    {\snorm*{#3}_{#4}}
    {
      \IfNoValueTF{#2}
        {\snorm{#3}_{#4}}
        {\snorm[#2]{#3}_{#4}}
    }
}

\NewDocumentCommand{\TNormIn}{s o m m}
{
  \IfBooleanTF{#1}
    {\tnorm*{#3}_{#4}}
    {
      \IfNoValueTF{#2}
        {\tnorm{#3}_{#4}}
        {\tnorm[#2]{#3}_{#4}}
    }
}

\NewDocumentCommand{\InnerIn}{s o m m m}
{
  \IfBooleanTF{#1}
    {\innerp*{#3}{#4}_{#5}}
    {
      \IfNoValueTF{#2}
        {\innerp{#3}{#4}_{#5}}
        {\innerp[#2]{#3}{#4}_{#5}}
    }
}

\NewDocumentCommand{\DualIn}{s o m m m}
{
  \IfBooleanTF{#1}
    {\dualp*{#3}{#4}_{#5}}
    {
      \IfNoValueTF{#2}
        {\dualp{#3}{#4}_{#5}}
        {\dualp[#2]{#3}{#4}_{#5}}
    }
}

\NewDocumentCommand{\NLp}{O{\defaultdomain} m m} { \NormIn{#2}{\Lp[#1]{#3}} }
\NewDocumentCommand{\NLpv}{O{\defaultdomain} m m} { \NormIn{#2}{\Lpv[#1]{#3}} }

\NewDocumentCommand{\NLtwo}{O{\defaultdomain} m} { \NormIn{#2}{\Ltwo[#1]} }
\NewDocumentCommand{\ILtwo}{O{\defaultdomain} m m} { \InnerIn{#2}{#3}{\Ltwo[#1]} }

\NewDocumentCommand{\NLtwov}{O{\defaultdomain} m} { \NormIn{#2}{\Ltwov[#1]} }
\NewDocumentCommand{\ILtwov}{O{\defaultdomain} m m} {\InnerIn{#2}{#3}{\Ltwov[#1]} }

\NewDocumentCommand{\NLinf}{O{\defaultdomain} m} { \NormIn{#2}{\Linf[#1]} }
\NewDocumentCommand{\NLinfv}{O{\defaultdomain} m} { \NormIn{#2}{\Linfv[#1]} }

\NewDocumentCommand{\NWmp}{O{\defaultdomain} m m m} { \NormIn{#2}{\Wmp[#1]{#3}{#4}} }
\NewDocumentCommand{\SNWmp}{O{\defaultdomain} m m m} { \SNormIn{#2}{\Wmp[#1]{#3}{#4}} }
\NewDocumentCommand{\TNWmp}{O{\defaultdomain} m m m} { \TNormIn{#2}{\Wmp[#1]{#3}{#4}} }

\NewDocumentCommand{\NWmpv}{O{\defaultdomain} m m m} { \NormIn{#2}{\Wmpv[#1]{#3}{#4}} }
\NewDocumentCommand{\SNWmpv}{O{\defaultdomain} m m m} { \SNormIn{#2}{\Wmpv[#1]{#3}{#4}} }
\NewDocumentCommand{\TNWmpv}{O{\defaultdomain} m m m} { \TNormIn{#2}{\Wmpv[#1]{#3}{#4}} }

\NewDocumentCommand{\NHm}{O{\defaultdomain} m m} { \NormIn{#2}{\Hm[#1]{#3}} }
\NewDocumentCommand{\SNHm}{O{\defaultdomain} m m} { \SNormIn{#2}{\Hm[#1]{#3}} }
\NewDocumentCommand{\TNHm}{O{\defaultdomain} m m} { \TNormIn{#2}{\Hm[#1]{#3}} }
\NewDocumentCommand{\IHm}{O{\defaultdomain} m m m} { \InnerIn{#2}{#3}{\Hm[#1]{#4}} }

\NewDocumentCommand{\NHmv}{O{\defaultdomain} m m} { \NormIn{#2}{\Hmv[#1]{#3}} }
\NewDocumentCommand{\SNHmv}{O{\defaultdomain} m m} { \SNormIn{#2}{\Hmv[#1]{#3}} }
\NewDocumentCommand{\TNHmv}{O{\defaultdomain} m m} { \TNormIn{#2}{\Hmv[#1]{#3}} }
\NewDocumentCommand{\IHmv}{O{\defaultdomain} m m m} { \InnerIn{#2}{#3}{\Hmv[#1]{#4}} }

\NewDocumentCommand{\NHone}{O{\defaultdomain} m} { \NormIn{#2}{\Hone[#1]} }
\NewDocumentCommand{\SNHone}{O{\defaultdomain} m} { \SNormIn{#2}{\Hone[#1]} }
\NewDocumentCommand{\TNHone}{O{\defaultdomain} m} { \TNormIn{#2}{\Hone[#1]} }
\NewDocumentCommand{\IHone}{O{\defaultdomain} m m} { \InnerIn{#2}{#3}{\Hone[#1]} }

\NewDocumentCommand{\NHonev}{O{\defaultdomain} m} { \NormIn{#2}{\Honev[#1]} }
\NewDocumentCommand{\SNHonev}{O{\defaultdomain} m} { \SNormIn{#2}{\Honev[#1]} }
\NewDocumentCommand{\TNHonev}{O{\defaultdomain} m} { \TNormIn{#2}{\Honev[#1]} }
\NewDocumentCommand{\IHonev}{O{\defaultdomain} m m} { \InnerIn{#2}{#3}{\Honev[#1]} }

\NewDocumentCommand{\NHdiv}{O{\defaultdomain} m} { \NormIn{#2}{\Hdiv[#1]} }
\NewDocumentCommand{\IHdiv}{O{\defaultdomain} m m} {\InnerIn{#2}{#3}{\Hdiv[#1]} }

\NewDocumentCommand{\NHcurl}{O{\defaultdomain} m} { \NormIn{#2}{\Hcurl[#1]} }
\NewDocumentCommand{\IHcurl}{O{\defaultdomain} m m} {\InnerIn{#2}{#3}{\Hcurl[#1]} }

\newcommand{\timeint}{(0,T)}
\newcommand{\ctimeint}{[0,T]}

\NewDocumentCommand{\BochnerLp}{O{\timeint} m m}
{
    L^{#3}(#1;#2)
}

\NewDocumentCommand{\BochnerHm}{O{\timeint} m m}
{
    H^{#3}(#1;#2)
}

\NewDocumentCommand{\BochnerC}{O{\ctimeint} m m}
{
    C^{#3}(#1;#2)
}

\NewDocumentCommand{\TLp}{O{\timeint} m m} { \BochnerLp[#1]{#2}{#3} }
\NewDocumentCommand{\TLtwo}{O{\timeint} m} { \BochnerLp[#1]{#2}{2} }
\NewDocumentCommand{\TLinf}{O{\timeint} m} {\BochnerLp[#1]{#2}{\infty} }
\NewDocumentCommand{\THm}{O{\timeint} m m} { \BochnerHm[#1]{#2}{#3} }
\NewDocumentCommand{\TCont}{O{\ctimeint} m} {\BochnerC[#1]{#2}{0} }

\NewDocumentCommand{\jump}{o m}
{
    \left\llbracket #2 \right\rrbracket
    \IfValueT{#1}{_{#1}}
}

\NewDocumentCommand{\avg}{o m}
{
    \left\{\!\!\left\{ #2 \right\}\!\!\right\}
    \IfValueT{#1}{_{#1}}
}

\newcommand{\rn}[1]{\uppercase\expandafter{\romannumeral #1}}

\newcommand{\bs}{\backslash}

\newcommand{\revise}[1]{{\color{blue}#1}}
